\documentclass[10pt,twoside,a4paper,reqno]{amsart}
\usepackage{amscd,amsmath,amsthm,amsfonts,latexsym,amssymb}

\theoremstyle{plain}
\newtheorem{theo+}           {Theorem}
\newtheorem{prop+}           {Proposition}
\newtheorem{coro+}           {Corollary}
\newtheorem{lemm+}           {Lemma}

\theoremstyle{definition}
\newtheorem{defi+}           {Definition}
\newtheorem{problem}         {Problem}

\theoremstyle{remark}
\newtheorem{rema+}           {Remark}

\newenvironment{theorem}{\begin{theo+}}{\end{theo+}}

\newenvironment{corollary}{\begin{coro+}}{\end{coro+}}
\newenvironment{lemma}{\begin{lemm+}}{\end{lemm+}}
\newenvironment{remark}{\begin{rema+}}{\end{rema+}}
\newenvironment{definition}{\begin{defi+}}{\end{defi+}}

\newcommand {\bC} {\mathbb {C}}
\newcommand {\bR} {\mathbb {R}}
\newcommand {\bZ} {\mathbb {Z}}
\newcommand {\bX} {\mathbb {X}}
\newcommand {\bN} {\mathbb {N}}
\newcommand {\bM} {\mathbb {M}}
\newcommand {\al} {\alpha}
\newcommand {\be} {\beta}

\newcommand {\de} {\delta}
\newcommand {\la} {\lambda}
\newcommand {\ze} {\zeta}
\newcommand {\si} {\sigma}
\newcommand {\om} {\omega}

\newcommand {\calC} {\mathcal {C}}

\newcommand {\calM} {\mathcal {M}}

\newcommand {\calZ} {\mathcal {Z}}
\newcommand {\calA} {\mathcal {A}}

\newcommand {\calP} {\mathcal {P}}
\newcommand{\co} {\text{conv}}
\newcommand{\tco} {\text{{\em conv}}}

\newcommand{\ba}{\mathbf a}
\newcommand{\bb}{\mathbf b}
\newcommand{\bc}{\mathbf c}

\newcommand{\pr}{\prec}
\newcommand{\dd}{\mathfrak{d}}
\newcommand{\1}{\mathbf{1}}
\newcommand{\su}{\text{supp}}
\newcommand{\tsu}{\text{{\em supp}}}

\begin{document}

\numberwithin{equation}{section}

\title[Choquet order, Lam\'e operators and orthogonal polynomials]
{Choquet order for spectra of higher Lam\'e operators and 
orthogonal polynomials}

\author[J.~Borcea]{Julius Borcea}
\address{Department of Mathematics, Stockholm University, SE-106 91 Stockholm,
   Sweden}
\email{julius@math.su.se}

\begin{abstract}
We establish a hierarchy of weighted majorization relations for 
the singularities of generalized Lam\'e equations and the zeros of their 
Van Vleck and Heine-Stieltjes polynomials as well as for 
multiparameter spectral polynomials of
higher Lam\'e operators. These relations translate into natural dilation 
and subordination properties in the Choquet order for certain probability 
measures associated with the aforementioned polynomials. As a 
consequence we obtain new inequalities for the moments and logarithmic 
potentials of the corresponding root-counting measures and their
weak-$^*$ limits in the semi-classical and various thermodynamic asymptotic 
regimes. We also prove analogous results for systems of orthogonal polynomials
such as Jacobi polynomials.
\end{abstract}

\maketitle

\section{Introduction}\label{s1}

The {\em generalized Lam\'e equation} in algebraic form is the second order
differential equation
\begin{equation}\label{2-lame} 
Q_2(z)y''(z)+Q_1(z)y'(z)+Q_0(z)y(z)=0,
\end{equation}
where $Q_2,Q_1,Q_0\in\bC[z]$ with $\deg Q_2=p$, $\deg Q_1=p-1$, 
$\deg Q_0\le p-2$. Particularly important cases are $p=2$ and $p=3$,
which correspond to the 
hypergeometric differential equation and Heun's equation, respectively 
(cf.~\cite{MS}). The classical 
Heine-Stieltjes multiparameter spectral problem deals with so-called
{\em Lam\'e solutions of the first kind} (of given degree and type) to
equation~\eqref{2-lame} and may be formulated as follows: given 
$Q_2(z), Q_1(z)$ as above and  
$n\in\bN$ find a polynomial $V(z)$ of degree at most $p-2$ and a polynomial 
$S(z)$
of degree $n$ such that~\eqref{2-lame} holds for $Q_0(z)=V(z)$ and 
$y(z)=S(z)$. If such
$V(z)$ and $S(z)$ exist we say that~\eqref{2-lame} is $n$-{\em solvable}. A 
generalized Lam\'e equation is {\em solvable} if it is $n$-solvable 
for all $n\in\bN$. The coefficients $V(z)$ are called 
{\em Van Vleck polynomials}
and the corresponding solutions $S(z)$ are known as 
{\em Heine-Stieltjes polynomials}. These two classes are also referred to 
as {\em Lam\'e polynomials} or {\em generalized spectral polynomials}
for~\eqref{2-lame}.

There are several known sufficient conditions for the solvability 
of equation~\eqref{2-lame}. For instance, Heine~\cite{H} proved that 
for any $n\in\bN$ there exist at most
$$\si(n):=\binom{n+p-2}{n}$$
different Van Vleck polynomials $V(z)$ for which~\eqref{2-lame} has a 
polynomial solution $y(z)=S(z)$ of degree $n$. Heine's text is written in a 
traditional XIXth century style German and the exact statements it contains 
seem to have created some confusion (cf.~\cite{Mar}). 
Szeg\"o~\cite[\S 6.8]{Sz} quotes
this result and adds that ``Heine asserts that, in general, there are exactly
$\si(n)$ determinations of this kind''. As explained in~\cite{BS,BBS}, Heine 
actually showed that if the coefficients of $Q_2(z)$ and $Q_1(z)$ are 
{\em algebraically independent} -- that is, these coefficients satisfy no 
algebraic equation
with integer coefficients -- then~\eqref{2-lame} is solvable. Moreover, if 
this is the case then for any $n\in\bN$ there exist exactly $\si(n)$ different
Van Vleck polynomials $V(z)$ of degree $p-2$ 
and the same number of corresponding monic
Heine-Stieltjes polynomials $S(z)$ of degree $n$. An explicit characterization
of the exceptional cases when this number is strictly less than $\si(n)$ 
seems to be lacking for the moment~\cite{Mar}. In the general case, 
Heine's arguments imply
that~\eqref{2-lame} is $n$-solvable for all sufficiently large $n$~\cite{BBS}. 

The solvability of~\eqref{2-lame} has been established under various other
assumptions, most notably when $Q_2(z)$ and $Q_1(z)$ have strictly interlacing
real zeros and the leading coefficient of $Q_1(z)$ is positive. This case 
is particularly interesting from a physical point of view and has attracted
a lot of attention in recent years. Indeed, differential equations of the 
form~\eqref{2-lame} whose coefficients satisfy the above condition arise 
naturally 
when separating variables in the Laplace equation in spherical coordinates and 
yield important examples of quantum completely integrable systems such as
generalized (real or complex) Gaudin spin chains~\cite{BJMT,B,BT}. 
A fundamental result of Stieltjes~\cite{St} -- also known as the 
Heine-Stieltjes theorem \cite{Sz} -- asserts that if $Q_2(z)$ and $Q_1(z)$ 
have strictly interlacing real zeros and $Q_1(z)$ has positive 
leading coefficient then for each $n\in\bN$ there are 
exactly $\si(n)$ different Van Vleck polynomials $V(z)$ of degree $p-2$ 
and the same number of corresponding monic 
Heine-Stieltjes polynomials $S(z)$ of 
degree $n$. The latter are given by all possible ways of 
distributing the zeros of $S(z)$ in the 
$p-1$ open intervals defined by the zeros of $Q_2(z)$. Stieltjes actually 
showed that the zeros of $S(z)$ are the coordinates
of the equilibrium points of a certain electrostatic potential. Similar 
results have recently been obtained in cases when $Q_2(z)$ has all real zeros
and the residues in the partial fractional decomposition of $Q_1(z)Q_2(z)^{-1}$
have mixed signs~\cite{DvA,G}. 

Let us assume that $Q_2(z)$ and $Q_1(z)$ are such that
\begin{equation}\label{not-1}
\begin{split}
&Q_2(z)=\prod_{l=1}^{p}(z-\ze_l)\text{ and }\frac{Q_1(z)}{Q_2(z)}
=\sum_{l=1}^{p}\frac{a_l}{z-\ze_l},\\
&\text{where }
\ze_l\in\bC\text{ and }a_l>0, 1\le l\le p.
\end{split}
\end{equation}
Note that if $Q_2(z)$ and $Q_1(z)$ are as above and equation~\eqref{2-lame} 
is solvable then any Van Vleck polynomial is of degree exactly $p-2$.  
P\'olya~\cite{P} -- and Klein and B\^ocher before him (cf.~\cite{DvA}) --
showed that in this case the zeros of all Van Vleck and 
Heine-Stieltjes polynomials lie in the convex hull of $\ze_1,\ldots,\ze_p$.
Extensions of this Gauss-Lucas type theorem to cases when the residues $a_i$
are not necessarily positive real numbers as well as various 
other results on the location of zeros of Lam\'e polynomials have since
been obtained~\cite{A,M,Z}. In this paper we show that much more is 
actually true. Namely, if~\eqref{not-1} holds then the 
zeros of any Van Vleck 
polynomial together with those of a corresponding Heine-Stieltjes polynomial 
and the zeros of $Q_2(z)$ satisfy 
a {\em weighted majorization relation} in the sense of~\cite{Bo} (see \S 2).
This amounts to a dilation property -- equivalently, a subordination relation
in the Choquet order -- for certain probability measures associated with the
generalized spectral polynomials and the singularities of 
equation~\eqref{2-lame}. A precise statement of this result 
is given in Theorem~\ref{t-1} below. As a consequence we obtain new 
inequalities for the moments and logarithmic potentials associated with 
the root-counting measures of Lam\'e polynomials and we establish similar
properties in the thermodynamic ($p\to\infty$) and semi-classical 
($n\to \infty$) asymptotic regimes (Corollaries~\ref{cor-1}--\ref{cor-3}). 
These results hold in the greatest possible generality and require
no additional assumptions besides~\eqref{not-1}. Therefore, Theorem~\ref{t-1}
and Corollary~\ref{cor-1} apply whenever equation~\eqref{2-lame} is 
$n$-solvable while 
Corollaries~\ref{cor-2}--\ref{cor-3} make sense in all cases 
when~\eqref{2-lame} is solvable and the considered limits exist (see 
\S \ref{s3} for several concrete examples). In the special case when 
$\ze_i\in \bR$, $1\le i\le p$, our results are a natural complement to those 
of~\cite{BS,BJMT,B,BT,MS} dealing with asymptotic distributions, limiting
level-spacings and mean densities of zeros of Lam\'e polynomials.

Various extensions of the Heine-Stieltjes multiparameter spectral problem to 
higher order linear ordinary differential operators with polynomial 
coefficients have been studied in~\cite{BS,BBS}. In particular, if $Q_2(z)$ 
and $Q_1(z)$ are as in~\eqref{not-1} and $k\ge 2$ then one may consider an
operator of the form
\begin{equation}\label{h-d}
\dd(z)=Q_2(z)\frac{d^k}{dz^k}+Q_1(z)\frac{d^{k-1}}{dz^{k-1}}.
\end{equation}
As in~\cite{BBS}, we call $\dd(z)$ a 
{\em higher order generalized Lam\'e operator} or a {\em higher Lam\'e 
operator} for short, provided that 
its {\em Fuchs index} $r:=p-k$ is non-negative. If $r=0$ then $\dd(z)$ is 
a so-called {\em hypergeometric type operator}.
Such operators 
and their polynomial eigenfunctions have important applications to the study 
of the Bochner-Krall problem and exactly solvable models 
(see, e.g., \cite{BR,BS,BBS} and references therein). The multiparameter 
spectral problem for a higher Lam\'e 
operator $\dd(z)$ is as follows: given $n\in\bN$ find a polynomial $V(z)$ of
degree at most $r$ such that the equation
\begin{equation}\label{k-lame}
\dd(z)y(z)+V(z)y(z)=0
\end{equation}
has a polynomial solution $y(z)=S(z)$ of degree $n$. One can then
define the notions of 
$n$-{\em solvability}, {\em solvability}, {\em higher Van Vleck} and 
{\em Heine-Stieltjes 
polynomials} -- that is, {\em higher spectral polynomials} or
{\em Lam\'e polynomials} -- corresponding to~\eqref{k-lame} by analogy with 
the 
terminology used for~\eqref{2-lame}. Several sufficient conditions 
for the solvability of~\eqref{k-lame} that extend those of Heine 
for~\eqref{2-lame} were recently obtained in ~\cite{BBS} (see \S 3). We show 
that 
whenever equation~\eqref{k-lame} is solvable its singularities and
the zeros of all corresponding higher Lam\'e polynomials satisfy weighted
majorization relations (Theorem~\ref{t-2}) and we establish natural analogs 
of Corollaries~\ref{cor-1}--\ref{cor-3} for the higher order case 
(Corollaries~\ref{cor-5}--\ref{cor-7}).

Our methods also yield interesting applications of the Choquet 
order$\,/\,$weighted majorization to the theory of 
orthogonal polynomials. In particular, we prove appropriate versions of the
aforementioned results for classical orthogonal polynomials such as 
(ultraspherical) Gegenbauer polynomials, (associated) Legendre polynomials, 
Chebyshev polynomials and indeed any family of Jacobi polynomials (see 
\S \ref{ss33}).

This paper is organized as follows. In \S \ref{s2} we recall the notion of 
weighted multivariate majorization from~\cite{Bo} as well as the 
definition and properties of the Choquet order for non-negative Radon 
measures. We state and prove our main results in \S \ref{s3}. In 
\S \ref{s4} we give several generalizations and discuss some related problems.

\section{Weighted Majorization and the Choquet Order}\label{s2}

The majorization preorder on $n$-tuples of
real numbers -- also known as the strong spectral order, vector majorization 
or classical majorization -- essentially 
quantifies the intuitive notion that the components of a real $n$-vector 
are less spread out than the components of another such vector. 
Several matrix versions of this notion have been proposed and studied in 
various contexts~\cite{MO}. 
A weighted multivariate extension of both 
classical and matrix majorization was introduced in~\cite{Bo}. In the special 
case of complex $n$-vectors the definition of {\em loc.~cit.~}is as follows.
For $m\in \bN$ set
\begin{equation}\label{sp-n}
\begin{split}
&\calA_{m}=\left\{\ba=(a_1,\ldots,a_m)\,\,\bigg|\,\,a_i\in [0,1],1\le i\le m,
\sum_{i=1}^{m}a_i=1\right\},\\
&\bX_{m}=\bC^{m}\times \calA_{m},\quad \bX=\bigcup_{n=1}^{\infty}\bX_{m}.
\end{split}
\end{equation}
Denote by $\co(\Omega)$ the (closed) convex hull of a (bounded) set 
$\Omega\subset \bC$ and by $X^T$ the transpose of a (row) vector 
$X=(x_1,\ldots,x_m)\in\bC^m$. We frequently write $\co(X)$ for
$\co(\{x_1,\ldots,x_m\})$. Let $\bM_{m,n}^{\text{rs}}$ be the set of all row 
stochastic $m\times n$ matrices. 

\begin{definition}\label{d-1}
The pair $(X,\ba)\in \bX_{m}$ is said to be {\em weightily majorized} by 
the pair $(Y,\bb)\in \bX_{n}$, denoted $(X,\ba)\pr (Y,\bb)$, if there  
exists a matrix $R\in \bM_{m,n}^{\text{rs}}$ such that 
$$\tilde{X}^T=R\tilde{Y}^T\text{ and }\,\bb=\ba R,$$
where $\tilde{X}^T$ and $\tilde{Y}^T$ are obtained by some 
(and then any) ordering of the coordinates of $X^T$ and $Y^T$, respectively.
\end{definition}

\begin{remark}\label{r-new1}
Note that if $(X,\ba)\pr (Y,\bb)$ then $X\in\co(Y)^{m}$ and 
the $\ba$-barycenter of 
$X$ must coincide with the $\bb$-barycenter of $Y$, that is, 
$\sum_{i=1}^{m}a_{i}x_i=\sum_{j=1}^{n}b_{j}y_j$. Moreover, it is 
clear that the weighted majorization relation 
is both reflexive and transitive, which makes it a preorder on $\bX$. One can
also show that for every $m\in\bN$ this preorder induces a partial order on  
the orbit space $\bC^m/\Sigma_m$, where $\Sigma_m$ is the symmetric group on
$m$ elements. 
\end{remark}

The following characterization of the weighted majorization relation  
may be found in~\cite[Theorem 1]{Bo}.

\begin{theorem}\label{t-m}
Let $(X,\ba)\in \bX_{m}$ and $(Y,\bb)\in \bX_{n}$, where 
$X=(x_1,\ldots,x_m)\in \bC^m$, $\ba=(a_1,\ldots,a_m)\in \calA_m$,
$Y=(y_1,\ldots,y_n)\in \bC^n$ and $\bb=(b_1,\ldots,b_n)\in \calA_n$.
The following conditions are equivalent:
\begin{itemize}
\item[(i)] for any (continuous) convex function $f:\bC\rightarrow \bR$ one 
has 
$$\sum_{i=1}^{m}a_{i}f(x_i)\le \sum_{j=1}^{n}b_{j}f(y_j);$$
\item[(ii)] the relation $(X,\ba)\pr (Y,\bb)$ holds.
\end{itemize}
\end{theorem}

\begin{remark}\label{r-1}
If $(X,\ba)\pr (Y,\bb)$ then the inequality in Theorem~\ref{t-m} (i)
holds for every convex function $f$ defined on $\co(Y)$.
\end{remark}

There is a natural connection between the weighted multivariate majorization 
relation and the Choquet order for non-negative Radon measures. The latter 
has been studied in the general context of locally convex separable 
topological vector spaces in e.g.~\cite{CFM} and subsequent papers. For 
measures defined on compact subsets of the complex plane the 
Choquet order and the main results of {\em op.~cit.~}may be described as 
follows.
Let $K$ be a convex compact subset of $\bC$, denote by $\calC(K)$ the space 
of real continuous functions on $K$ and let $\calP(K)$ the 
subset of $\calC(K)$ consisting of convex functions. If $\mu$ is a 
non-negative Radon measure on
$K$ and $f$ is a function on $K$ one defines $\mu(f)=\int_K f(y)d\mu(y)$.
The mass of $\mu$ is therefore $\mu(1)=\int_K d\mu(y)$ and if $\mu(1)>0$
then the barycenter of $\mu$ is the point 
$r(\mu):=\mu(1)^{-1}\int_K yd\mu(y)$. 

\begin{definition}\label{d-cho}
Given two non-negative Radon measures $\mu$ and $\nu$ on $K$ one says that
$\nu$ {\em dominates} $\mu$ {\em in the Choquet order} or that $\nu$ is a 
{\em dilation} of $\mu$, denoted $\mu\pr \nu$, if 
$\mu(f)\le \nu(f)$ for any $f\in\calP(K)$.
\end{definition}

The use of the term ``dilation'' in Definition~\ref{d-cho} is motivated by 
Definition~\ref{d-dil} and Theorem~\ref{t-dil} (iii) below. To formulate this 
result we need a few more concepts and 
notations. Let $\calM(K)$ be the set of all probability measures with 
$\su\mu\subseteq K$. Note that if $\mu\in\calM(K)$ then its barycenter 
$r(\mu)$ lies in $K$. 

\begin{remark}\label{r-2}
As is well known, the set $\calM(K)$ equipped with the weak-$^*$ topology is a
sequentially compact Hausdorff space. This will allow us to choose a 
convergent subsequence from any sequence of measures belonging to $\calM(K)$. 
\end{remark}

\begin{definition}\label{d-dil}
A {\em dilation} on $K$ is a weakly Borel measurable $\calM(K)$-valued 
function on $K$ that inverts the barycenter mapping. In other words, 
a map $T:K\to\calM(K)$, $x\mapsto T_x$, is a dilation on $K$ if $r(T_x)=x$ for
all $x\in K$ and the real-valued function on $K$ given by $x\mapsto T_x(f)$ 
is borelian for any $f\in\calC(K)$.
\end{definition}

If $T$ is a dilation on $K$ then for any non-negative Radon measure $\mu$ on
$K$ one can define a new such measure $\nu:=T(\mu)$ by setting
\begin{equation}\label{T-eq}
\nu(f)=\int_K T_x(f)d\mu(x),\quad f\in \calC(K).
\end{equation}
It is not difficult to show that the real-valued function on $K$ given by 
$x\mapsto T_x(f)$ is borelian and bounded whenever $f$ is a bounded borelian 
real-valued function on $K$ and that~\eqref{T-eq} actually 
holds for all such functions (cf.~\cite{CFM}). The main results of 
{\em loc.~cit.~}provide 
various descriptions of the Choquet order in a quite general setting. In the 
case discussed above these may be summarized as follows (see also \cite{EH}).

\begin{theorem}\label{t-dil}
If $K$ is a convex compact subset of $\bC$ and $\mu,\nu\in\calM(K)$ then  
the following conditions are equivalent:
\begin{itemize}
\item[(i)] $\mu\pr \nu$; 
\item[(ii)] for every convex combination $\mu=\sum_{i=1}^{n}\la_i\mu_i$
with $\mu_i\in\calM(K)$, $1\le i\le n$, there exists a corresponding convex 
combination $\nu=\sum_{i=1}^{n}\la_i\nu_i$ such that $\nu_i\in\calM(K)$
and $r(\nu_i)=r(\mu_i)$, $1\le i\le n$;
\item[(iii)] $\nu=T(\mu)$, where $T$ is a dilation on $K$.
\end{itemize}
\end{theorem}

\begin{remark}\label{r-new2}
If the conditions in Theorem~\ref{t-dil} hold then 
$\su(\mu)\subseteq \co(\su(\nu))$. This may be viewed as a ``continuous''
version of the corresponding result for weighted majorization 
(cf.~Remark~\ref{r-new1}).
\end{remark}

\section{Main Results and Proofs}\label{s3}

Given a complex polynomial $P$ of degree $d\ge 1$ we let $\calZ(P)$ be 
the $d$-tuple (or multiset) consisting of the zeros of $P$, where it is 
understood that each zero occurs as many times
as its multiplicity. In particular, $|\calZ(P)|=d$. To $P$ we associate
its root-counting measure, namely the (finite) real probability measure
given by
$$\mu_{_P}=|\calZ(P)|^{-1}\sum_{\ze\in\calZ(P)}\de_{\ze},$$
where $\de_{\ze}$ is the Dirac measure supported at $\ze$. The symbol $\vee$ 
is used below for the concatenation operation, that is,
if $(x_1,\ldots,x_m)\in\bC^m$ and $(y_1,\ldots,y_n)\in \bC^n$ then
\begin{equation*}
(x_1,\ldots,x_m)\vee (y_1,\ldots,y_n)=(x_1,\ldots,x_m,y_1,\ldots,y_n)
\in \bC^{m+n},
\end{equation*}
and the ``all ones'' vector is denoted by $\1_m=(1,\ldots,1)\in\bR^m$.

\subsection{Generalized Lam\'e operators}\label{ss31}

Suppose that $n\ge 2$ is an integer such that~\eqref{2-lame} is $n$-solvable 
and that 
$Q_2(z)$, $Q_1(z)$ satisfy~\eqref{not-1}. Let $S(z)$ be a Heine-Stieltjes 
polynomial of degree $n$ corresponding to a Van Vleck polynomial $V(z)$, so
that $\deg V=p-2$ (cf.~\S \ref{s1}). Let $\al=\al(n,p):=n-1+\sum_{l=1}^{p}a_l$ 
and define the following weight vectors:
\begin{equation}\label{a-2}
\begin{split}
&\ba=\frac{\al}{(p-1)\al+n-1}\1_{p-2},\quad 
\bb=\frac{\al+n-1}{n[(p-1)\al+n-1]}\1_n,\\
&\bc=\left(\frac{\al-a_1}{(p-1)\al+n-1},
\ldots,\frac{\al-a_p}{(p-1)\al+n-1}\right).\\
\end{split}
\end{equation}
Recall~\eqref{not-1} and note that $\al>1$, $\bc\in\calA_p$ while 
$\ba\vee\bb\in \calA_{n+p-2}$. Finally, set
\begin{equation}\label{z-2}
\calZ(V)=(v_1,\ldots,v_{p-2}),\quad \calZ(S)=(s_1,\ldots,s_n).
\end{equation}

We can now state our first main result.

\begin{theorem}\label{t-1}
With the above notations and assumptions the inequality
\begin{equation}\label{res-2}
\sum_{i=1}^{p-2}f(v_i)+
\left[1-\left(1-\frac{1}{n}\right)\!\left(1-\frac{1}{\al}\right)\right]
\sum_{j=1}^{n}f(s_j)\le 
\sum_{l=1}^{p}\left(1-\frac{a_l}{\al}\right)f(\ze_l)
\end{equation}
holds for any convex function $f:\bC\to \bR$ and if equality occurs 
in~\eqref{res-2} for some strictly convex function $f$ then the zeros
of $Q_2$ must be collinear. Equivalently, 
\begin{equation*}
\big(\calZ(V)\vee \calZ(S),\ba\vee\bb\big)\pr \big(\calZ(Q_2),\bc\big).
\end{equation*}
Thus there exists a matrix $R=R(n,p)\in\bM_{n+p-2,p}^{\text{rs}}$ such that
\begin{equation*}
\big(\calZ(V)\vee \calZ(S)\big)^T=R\calZ(Q_2)^T\text{ and }\bc=(\ba\vee\bb)R.
\end{equation*}
\end{theorem} 

\begin{remark}\label{r-appl-2}
As pointed out in \S \ref{s1}, the only requirements for Theorem~\ref{t-1} are
that \eqref{not-1} holds and equation \eqref{2-lame} is 
$n$-solvable. For instance, Stieltjes' theorem shows 
that \eqref{2-lame} is always solvable if 
$\calZ(Q_2)\subset \bR$ while Heine's result \cite{H} asserts that the same 
is true whenever $Q_2(z)$ and $Q_1(z)$ are algebraically independent.
\end{remark}

Note that in particular Theorem~\ref{t-1} immediately implies the 
P\'olya-Klein-B\^ocher result mentioned in \S \ref{s1}, namely 
$\calZ(V)\cup\calZ(S)\subseteq \co(\calZ(Q_2))$ (cf.~Remark~\ref{r-new1}). 
Now given a compact set $K\subset \bC$ and $\mu\in\calM(K)$ let
$$\mu^{(m)}:=\int|w|^md\mu(w),\quad m\in\bZ_+,$$
denote the moments of $\mu$. As is well known, the logarithmic 
potential of $\mu$ 
$$U^{\mu}(z)=\int\!\log|z-w|d\mu(w)$$
is subharmonic in $\bC$ and $U^{\mu}(z)=-\infty$ for every atom $z$ of $\mu$.

Clearly, \eqref{res-2} may be reformulated in terms of the Choquet
order for atomic probability measures with finite point spectrum: 

\begin{corollary}\label{cor-1}
In the situation of the preceding theorem one has
$$\frac{(p-2)\al}{(p-1)\al+n-1}\mu_{_V}+\frac{\al+n-1}{(p-1)\al+n-1}\mu_{_S}
\pr \tilde{\mu}_{_{Q_2}},$$
where $\mu_{_V}$ and $\mu_{_S}$ are the root-counting measures of $V$ and $S$,
respectively, while $\tilde{\mu}_{_{Q_2}}\in\calM(\tco(\calZ(Q_2)))$ is 
defined by 
$$\tsu\big(\tilde{\mu}_{_{Q_2}}\big)=\calZ(Q_2)=\{\ze_l\}_{l=1}^{p}\text{ and }
\tilde{\mu}_{_{Q_2}}(\{\ze_l\})=\frac{\al-a_l}{(p-1)\al+n-1}$$
for $1\le l\le p$. In particular,
$$\frac{(p-2)\al}{(p-1)\al+n-1}\mu_{_V}^{(m)}
+\frac{\al+n-1}{(p-1)\al+n-1}\mu_{_S}^{(m)}
\le \tilde{\mu}_{_{Q_2}}^{(m)}$$
for all $m\in\bZ_+$ and 
$$\frac{(p-2)\al}{(p-1)\al+n-1}U^{\mu_{_V}}(z)
+\frac{\al+n-1}{(p-1)\al+n-1}U^{\mu_{_S}}(z)\ge
U^{\tilde{\mu}_{_{Q_2}}}(z)$$
whenever $z\in \bC\setminus \tco(\calZ(Q_2))$.
\end{corollary}

In the semi-classical asymptotic regime ($n\to\infty$) 
Theorem~\ref{t-1} yields:

\begin{corollary}\label{cor-2}
Assume that \eqref{not-1} holds and that \eqref{2-lame} is solvable. Let 
$\{S_n(z)\}_{n\in\bN}$ be a sequence of monic Heine-Stieltjes polynomials
such that $\deg S_n=n$, $n\in\bN$, and $\{V_n(z)\}_{n\in\bN}$ be a 
corresponding sequence of Van Vleck polynomials with $\deg V_n=p-2$ 
normalized so that each $V_n$ is monic. Then
$$\frac{p-2}{p}\mu_{_V}^{*}+\frac{2}{p}\mu_{_S}^{*}\pr \tilde{\mu}_{_{Q_2}},$$
where $\mu_{_V}^{*}=\lim_{\Lambda\ni n\to\infty}^{*}\mu_{_{V_n}}$ and
$\mu_{_S}^{*}=\lim_{\Lambda\ni n\to\infty}^{*}\mu_{_{S_n}}$ for an  
appropriately chosen $\Lambda\subset \bN$. Equivalently,
there exists a dilation $T$ on $\tco(\calZ(Q_2))$ such that 
$p\tilde{\mu}_{_{Q_2}}=T\big((p-2)\mu_{_V}^{*}+2\mu_{_S}^{*}\big)$. 
In particular,
$$\frac{p-2}{p}\mu_{_V}^{*^{(m)}}
+\frac{2}{p}\mu_{_S}^{*^{(m)}}\le \tilde{\mu}_{_{Q_2}}^{(m)},
\quad m\in\bZ_+,$$ 
and 
$$\frac{p-2}{p}U^{\mu_{_V}^{*}}(z)
+\frac{2}{p}U^{\mu_{_S}^{*}}(z)\ge U^{\tilde{\mu}_{_{Q_2}}}(z)$$
for any $z\in \bC\setminus \tco(\calZ(Q_2))$.
\end{corollary}

Finally, we may also let $p\to\infty$ and consider various so-called 
thermodynamic asymptotic regimes (cf., e.g., \cite{BJMT,B,BT}). In this
case we get the following result.

\begin{corollary}\label{cor-3}
Let $\{Q_{2,p}(z)\}_{p=2}^{\infty}$ and $\{Q_{1,p}(z)\}_{p=2}^{\infty}$ be
two sequences of polynomials such that for all $p\ge 2$ the pair 
$(Q_{2,p}(z),Q_{1,p}(z))$ satisfies \eqref{not-1} and the corresponding 
equation \eqref{2-lame} is solvable. Assume further that $K\subset \bC$ is a 
compact set with  
$\calZ(Q_{2,p})\subset K$, $p\ge 2$, and that $\{S_{p,n}(z)\}_{n\in\bN}$,
respectively $\{V_{p,n}(z)\}_{n\in\bN}$, is a sequence of monic 
Heine-Stieltjes polynomials, respectively normalized Van Vleck polynomials,
associated with the resulting system of equations such that 
$\deg S_{p,n}=n$, $\deg V_{p,n}=p-2$ and $V_{p,n}$ is monic for all
$p\ge 2$ and $n\in\bN$. Then
$$_{*}\mu_{_V}\pr {_{*}\mu}_{_{Q_2}},$$  
where $_{*}\mu_{_V}=\lim_{\Gamma\ni p\to \infty}^{*}\mu_{_{V_{p,n(p)}}}$ 
and ${_{*}\mu}_{_{Q_2}}=\lim_{\Gamma\ni p\to \infty}^{*}\mu_{_{Q_{2,p}}}$ for 
an appropriately chosen $\Gamma\subset \bN$. Hence there exists
a dilation $T$ on $K$ such that ${_{*}\mu}_{_{Q_2}}=T\big({_{*}\mu}_{_V}\big)$.
In particular, 
$${_{*}\mu}_{_V}^{(m)}\le {_{*}\mu}_{_{Q_2}}^{(m)},\quad m\in\bZ_+,$$ 
while 
$$U^{{_{*}\mu}_{_V}}(z)\ge U^{{_{*}\mu}_{_{Q_2}}}(z)$$ 
whenever $z\in \bC\setminus K$.
\end{corollary}

\subsection{Higher Lam\'e operators}\label{ss32}

Let now $k\ge 2$ be a fixed integer and consider an order $k$ generalized 
Lam\'e operator $\dd(z)$ with Fuchs index $r:=p-k$ as in \eqref{h-d} and the
corresponding multiparameter spectral problem~\eqref{k-lame}. Assume that the
latter is $n$-solvable and that $(S(z),V(z))$ is a pair of (higher) spectral 
polynomials with $\deg S=n$ and $\deg V=r$ satisfying \eqref{k-lame}.
Let $\al_k=\al(n,p,k):=n-k+1+\sum_{l=1}^{p}a_l$ 
and define the following weight vectors:
\begin{equation}\label{a-k}
\begin{split}
&\ba_k=\frac{\al_k}{(p-1)\al_k+n-k+1}\1_{r},\quad 
\bb_k=\frac{(k-1)\al_k+n-k+1}{n[(p-1)\al_k+n-k+1]}\1_n,\\
&\bc_k=\left(\frac{\al_k-a_1}{(p-1)\al_k+n-k+1},
\ldots,\frac{\al_k-a_p}{(p-1)\al_k+n-k+1}\right).\\
\end{split}
\end{equation}
One clearly has $\al_k>1$, $\bc_k\in \calA_p$, $\ba_k\vee \bb_k\in\calA_{n+r}$,
$\ba_2=\ba$, $\bb_2=\bb$ and $\bc_2=\bc$. Since in this case $\deg V=r$ we 
adapt notation \eqref{z-2} to the current situation simply by setting 
$\calZ(V)=(v_1,\ldots,v_{r})$. 

The analog of Theorem~\ref{t-1} for higher Lam\'e operators reads as follows.

\begin{theorem}\label{t-2}
Under the above assumptions the inequality
\begin{equation}\label{res-k}
\sum_{i=1}^{r}f(v_i)+
\left[1-\left(1-\frac{k-1}{n}\right)\!\left(1-\frac{1}{\al_k}\right)\right]
\sum_{j=1}^{n}f(s_j)\le 
\sum_{l=1}^{p}\left(1-\frac{a_l}{\al_k}\right)f(\ze_l)
\end{equation}
holds for any convex function $f:\bC\to \bR$ and if equality occurs 
in~\eqref{res-k} for some strictly convex function $f$ then the $\ze_l$'s
must be collinear. Equivalently, 
\begin{equation*}
\big(\calZ(V)\vee \calZ(S),\ba_k\vee\bb_k\big)\pr \big(\calZ(Q_2),\bc_k\big).
\end{equation*}
Thus there exists a matrix $R_k=R(n,p,k)\in\bM_{n+r,p}^{\text{rs}}$ such that
\begin{equation*}
\big(\calZ(V)\vee \calZ(S)\big)^T=R_k\calZ(Q_2)^T
\text{ and }\bc_k=(\ba_k\vee\bb_k)R_k.
\end{equation*}
\end{theorem}

\begin{remark}\label{r-strong}
The proof of Theorem~\ref{t-2} actually yields an inequality stronger
than~\eqref{res-k} involving all four zero sets $\calZ(V)$, $\calZ(S)$,
$\calZ(S^{(k-1)})$ and $\calZ(Q_2)$ (see~\eqref{eq-strong} below). 
\end{remark}

\begin{remark}\label{r-appl-k}
Theorem~\ref{t-2} applies to all situations when $Q_2(z),Q_1(z)$ satisfy 
\eqref{not-1} and equation \eqref{k-lame} is $n$-solvable. By 
\cite[Theorem 5]{BBS} this is always true for all sufficiently large $n$.
Moreover, it was shown in {\em op.~cit.~}that \eqref{k-lame} is
$n$-solvable for any $n\in\bN$ in each of the following cases: 
(i) $Q_2(z)$ and $Q_1(z)$ are algebraically independent (ii) $\calZ(Q_2)\subset
\bR$ (iii) $\dd(z)$ is a hyperbolicity preserving 
operator (HPO for short), i.e., it maps polynomials with all real zeros to 
polynomials with all real zeros. A complete classification of all HPOs was 
recently obtained in \cite{BBS1} (see also \cite{BBS2}).
\end{remark}

Natural extensions of Corollaries~\ref{cor-1}--\ref{cor-3} to higher Lam\'e
operators are as follows.

\begin{corollary}\label{cor-5}
In the situation of Theorem~\ref{t-2} one has
$$\frac{(p-k)\al_k}{(p-1)\al_k+n-1}\mu_{_V}
+\frac{(k-1)\al_k+n-1}{(p-1)\al_k+n-1}\mu_{_S}
\pr \tilde{\mu}_{_{Q_2}},$$
where $\mu_{_V}$ and $\mu_{_S}$ are the root-counting measures of $V$ and $S$,
respectively, while $\tilde{\mu}_{_{Q_2}}\in\calM(\tco(\calZ(Q_2)))$ is 
defined by 
$$\tsu\big(\tilde{\mu}_{_{Q_2}}\big)=\calZ(Q_2)=\{\ze_l\}_{l=1}^{p}\text{ and }
\tilde{\mu}_{_{Q_2}}(\{\ze_l\})=\frac{\al_k-a_l}{(p-1)\al_k+n-1}$$
for $1\le l\le p$. In particular,
$$\frac{(p-k)\al_k}{(p-1)\al_k+n-1}\mu_{_V}^{(m)}
+\frac{(k-1)\al_k+n-1}{(p-1)\al_k+n-1}\mu_{_S}^{(m)}
\le \tilde{\mu}_{_{Q_2}}^{(m)}$$
for all $m\in\bZ_+$ and 
$$\frac{(p-k)\al_k}{(p-1)\al_k+n-1}U^{\mu_{_V}}(z)
+\frac{(k-1)\al_k+n-1}{(p-1)\al_k+n-1}U^{\mu_{_S}}(z)\ge
U^{\tilde{\mu}_{_{Q_2}}}(z)$$
whenever $z\in \bC\setminus \tco(\calZ(Q_2))$.
\end{corollary}

\begin{corollary}\label{cor-6}
Assume that \eqref{not-1} holds and that \eqref{k-lame} is solvable. Let 
$\{S_n(z)\}_{n\in\bN}$ be a sequence of monic higher Heine-Stieltjes 
polynomials
such that $\deg S_n=n$, $n\in\bN$, and $\{V_n(z)\}_{n\in\bN}$ be a 
corresponding sequence of higher Van Vleck polynomials with $\deg V_n=p-2$ 
normalized so that each $V_n$ is monic. Then
$$\frac{p-k}{p}\mu_{_V}^{*}+\frac{k}{p}\mu_{_S}^{*}\pr \tilde{\mu}_{_{Q_2}},$$
where $\mu_{_V}^{*}=\lim_{\Lambda\ni n\to\infty}^{*}\mu_{_{V_n}}$ and
$\mu_{_S}^{*}=\lim_{\Lambda\ni n\to\infty}^{*}\mu_{_{S_n}}$ for an 
appropriately chosen $\Lambda\subset \bN$. Equivalently,
there exists a dilation $T$ on $\tco(\calZ(Q_2))$ such that 
$p\tilde{\mu}_{_{Q_2}}=T\big((p-k)\mu_{_V}^{*}+k\mu_{_S}^{*}\big)$. 
In particular, 
$$\frac{p-k}{p}\mu_{_V}^{*^{(m)}}
+\frac{k}{p}\mu_{_S}^{*^{(m)}}\le \tilde{\mu}_{_{Q_2}}^{(m)},
\quad m\in\bZ_+,$$ 
and 
$$\frac{p-k}{p}U^{\mu_{_V}^{*}}(z)
+\frac{k}{p}U^{\mu_{_S}^{*}}(z)\ge U^{\tilde{\mu}_{_{Q_2}}}(z)$$
for any $z\in \bC\setminus \tco(\calZ(Q_2))$.
\end{corollary}

\begin{corollary}\label{cor-7}
Let $\{Q_{2,p}(z)\}_{p=2}^{\infty}$ and $\{Q_{1,p}(z)\}_{p=2}^{\infty}$ be
two sequences of polynomials such that for all $p\ge 2$ the pair 
$(Q_{2,p}(z),Q_{1,p}(z))$ satisfies \eqref{not-1} and the corresponding 
higher Lam\'e equation \eqref{k-lame} is solvable. Assume further that 
$K\subset \bC$ is a compact set with  
$\calZ(Q_{2,p})\subset K$, $p\ge 2$, and that $\{S_{p,n}(z)\}_{n\in\bN}$,
respectively $\{V_{p,n}(z)\}_{n\in\bN}$, is a sequence of monic higher  
Heine-Stieltjes polynomials, respectively normalized higher 
Van Vleck polynomials,
associated with the resulting system of higher Lam\'e equations such that 
$\deg S_{p,n}=n$, $\deg V_{p,n}=p-2$ and $V_{p,n}$ is monic for all
$p\ge 2$ and $n\in\bN$. Then
$$_{*}\mu_{_V}\pr {_{*}\mu}_{_{Q_2}},$$ 
where $_{*}\mu_{_V}=\lim_{\Gamma\ni p\to \infty}^{*}\mu_{_{V_{p,n(p)}}}$ 
and ${_{*}\mu}_{_{Q_2}}=\lim_{\Gamma\ni p\to \infty}^{*}\mu_{_{Q_{2,p}}}$ for 
an appropriately chosen $\Gamma\subset \bN$. Hence there exists
a dilation $T$ on $K$ such that ${_{*}\mu}_{_{Q_2}}=T\big({_{*}\mu}_{_V}\big)$.
In particular, 
$${_{*}\mu}_{_V}^{(m)}\le {_{*}\mu}_{_{Q_2}}^{(m)},\quad m\in\bZ_+,$$ 
while 
$$U^{{_{*}\mu}_{_V}}(z)\ge U^{{_{*}\mu}_{_{Q_2}}}(z)$$ 
whenever $z\in \bC\setminus K$.
\end{corollary}

\subsection{Orthogonal polynomials}\label{ss33}

As we shall now explain, the above results have interesting yet apparently 
unknown analogs for important classes of orthogonal polynomials such as
Jacobi polynomials. Recall that the latter are defined by
$$P_n^{(\al,\be)}(z)=\frac{(-1)^n}{2^n n!}(1-z)^{-\al}(1+z)^{-\be}
\frac{d^n}{dz^n}\left[(1-z)^{n+\al}(1+z)^{n+\be}\right],$$
where $\al>-1$, $\be>-1$. For special values of the parameters $\al$ and 
$\be$ one 
gets (up to a normalizing factor) all the other classical Jacobi-like 
polynomials
including (ultraspherical) Gegenbauer polynomials, (associated) Legendre
polynomials and Chebyshev polynomials of the first or second kind. The 
relative location and asymptotic behaviour of the zeros of 
Jacobi polynomials have been of permanent interest in view of their important
role as nodes of Gaussian quadrature formulae and their nice electrostatic
interpretation (Bethe Ansatz) \cite{Sz}. 

We prove the following result.

\begin{theorem}\label{t-j}
Let $\calZ\big(P_n^{(\al,\be)}\big)=\{\ze_{n,i}\}_{i=1}^{n}$ be the zero
set of $P_n^{(\al,\be)}$. Then 
$$\frac{1}{n}\sum_{i=1}^{n}f(\ze_{n,i})\le \frac{(n+\be)f(1)
+(n+\al)f(-1)}{2n+\al+\be}$$
for any convex function $f:[-1,1]\to\bR$.
\end{theorem}

\begin{remark}\label{r-j-strong}
The proof of Theorem~\ref{t-j} yields in fact an even stronger relation that 
involves both zero sets $\calZ\big(P_n^{(\al,\be)}\big)$ and 
$\calZ\Big(P_n^{(\al,\be)'}\Big)$, see \eqref{eq-j-strong} in \S \ref{ss34}
below.
\end{remark}

\begin{remark}\label{r-jac}
Note that if $\mu_{_{P_n^{(\al,\be)}}}$ denotes the root-counting measure of 
$P_n^{(\al,\be)}$ then Theorem~\ref{t-j} may be restated in terms of the 
Choquet order simply as
$$\mu_{_{P_n^{(\al,\be)}}}\pr\frac{n+\be}{2n+\al+\be}\de_1+
\frac{n+\al}{2n+\al+\be}\de_{-1}.$$
In particular, by letting $n\to\infty$ we get
\begin{equation}\label{eq-as-j}
\mu_{_{\al,\be}}^{*}\pr\frac{1}{2}(\de_1+\de_{-1}),
\end{equation}
where $\mu_{_{\al,\be}}^{*}:=\lim_{n\to\infty}^{*}\mu_{_{P_n^{(\al,\be)}}}$
is the $^*$-limiting distribution of the zeros of $P_n^{(\al,\be)}$. As
is well known (see, e.g., \cite{Sz}) the latter is the (uniform) arcsine 
distribution and thus \eqref{eq-as-j} may be rewritten as
$$\frac{1}{\pi}\int_{-1}^{1}\frac{f(z)}{\sqrt{1-z^2}}dz\le 
\frac{f(1)+f(-1)}{2}$$
for any convex function $f:[-1,1]\to\bR$. Although elementary, the above 
inequality is not completely obvious; arguably the most direct way of proving
it is to note that it actually amounts to showing that
$$\frac{1}{\pi}\int_{-1}^{1}\frac{|z-c|}{\sqrt{1-z^2}}dz\le \max(1,|c|)$$
for any $c\in\bR$, which is a trivial exercise.
\end{remark}

\subsection{Proofs}\label{ss34}

Fix an integer $m\ge 2$ and let $z_i$, $1\le i\le m$, be (not necessarily
distinct) points in the complex plane that do not coalesce into a single one.
Given $\tau_i>0$, $1\le i\le m$, such that $\sum_{i=1}^{m}\tau_i=1$ we define a
meromorphic function
\begin{equation}\label{phi}
\Phi(z)=\sum_{i=1}^{m}\frac{\tau_i}{z-z_i}.
\end{equation}
Functions of this type are sometimes called generalized derivatives in the
sense of Sz.-Nagy (see \cite{Bo}) and may be interpreted as the resulting 
electrostatic force of a planar charge configuration (cf.~{\em loc.~cit.}). 
One of the key ingredients
in our proofs is \cite[Theorem 2]{Bo}, which we restate as follows:

\begin{lemma}\label{l-1}
Let $\Phi$ be as in \eqref{phi} and denote its zeros by $w_j$, 
$1\le j\le m-1$, where it is understood that $z_i$ counts as a ``zero'' of 
$\Phi$ of multiplicity $m_i-1$ if it occurs precisely 
$m_i$ times in \eqref{phi}. Then
\begin{equation}\label{eq-phi}
\sum_{j=1}^{m-1}f(w_j)\le \sum_{i=1}^{m}(1-\tau_i)f(z_i)
\end{equation}
for any convex function $f:\bC\to \bR$.
\end{lemma}

\begin{remark}\label{r-eq}
The arguments in \cite{Bo} further imply that if $f$ is a strictly convex
function such that equality is attained in~\eqref{eq-phi} then the $z_i$'s
must be collinear.
\end{remark}

We emphasize an important special case of Lemma~\ref{l-1}:

\begin{corollary}\label{c-cons}
If $P\in\bC[z]$ is such that $\deg P=d\ge 2$ then
$$(d-i+1)\sum_{w\in\calZ(P^{(i)})}f(w)\le 
(d-i)\sum_{z\in\calZ(P^{(i-1)})}f(z)$$
for any convex function $f:\bC\to \bR$ and $1\le i\le d-1$.
\end{corollary}

To prove Theorem~\ref{t-2} let $(S(z),V(z))$ be a pair of (higher) spectral 
polynomials of $\dd(z)$ as in \S \ref{ss32}, set
$$\calZ\big(S^{(i)}\big)
=\big(s_1^{(i)},\ldots,s_{n-i}^{(i)}\big),\,1\le i\le n-1,$$
and note that by \eqref{not-1} equation \eqref{k-lame} may be rewritten as
\begin{equation}\label{gen-k}
-\frac{V(z)S(z)}{Q_2(z)S^{(k-1)}(z)}=\frac{S^{(k)}(z)}{S^{(k-1)}(z)}+
\frac{Q_1(z)}{Q_2(z)}=\sum_{j=1}^{n-k+1}\frac{1}{z-s_j^{(k-1)}}+
\sum_{l=1}^{p}\frac{a_l}{z-\ze_l}.
\end{equation}
Since $\calZ(VS)=\calZ(V)\vee \calZ(S)$ (by the convention made at the 
beginning of \S \ref{s3}) and $a_l>0$, $1\le l\le p$, we deduce
from Lemma~\ref{l-1} that for any convex function $f:\bC\to \bR$ one has
\begin{equation}\label{eq-strong}
\sum_{i=1}^{r}f(v_i)+\sum_{j=1}^{n}f(s_j)\le 
\left(1-\frac{1}{\al_k}\right)\sum_{j=1}^{n-k+1}f\!\left(s_j^{(k-1)}\!\right)
+\sum_{l=1}^{p}\left(1-\frac{a_l}{\al_k}\right)f(\ze_l).
\end{equation}
On the other hand, by Corollary~\ref{c-cons} we know that 
$$n\sum_{j=1}^{n-k+1}f\!\left(s_j^{(k-1)}\!\right)\le 
(n-k+1)\sum_{j=1}^{n}f(s_j),$$
which combined with~\eqref{eq-strong} yields~\eqref{res-k} after some 
straightforward computations. Theorem~\ref{t-1} follows from the above simply 
by letting $k=2$. 

Since functions of the type
$\bC\ni w\mapsto |w|^m$, $m\in\bZ_+$, and
$\co(\calZ(Q_2))\ni w\mapsto -\log|z-w|$, 
$z\notin \co(\calZ(Q_2))$, are convex,
Corollaries~\ref{cor-1} and~\ref{cor-5} are immediate
consequences of Theorems~\ref{t-1} and~\ref{t-2}, respectively. 
Corollaries~\ref{cor-2} and~\ref{cor-6} follow from Corollaries~\ref{cor-1}
and~\ref{cor-5}, respectively, by using Remark~\ref{r-2}
and the expression for $\al_k$ (i.e., 
$\al_k=\al(n,p,k)=n-k+1+\sum_{l=1}^{p}a_l$) and noticing that
$$\frac{(p-k)\al_k(n,p,k)}{(p-1)\al_k(n,p,k)+n-1}\to\frac{p-k}{p}\,
\text{ and }\,
\frac{(k-1)\al_k(n,p,k)+n-1}{(p-1)\al_k(n,p,k)+n-1}\to \frac{k}{p}$$
as $n\to \infty$, $p$ being fixed. 
To prove Corollary~\ref{cor-7} (hence also Corollary~\ref{cor-3}) 
we use again Corollary~\ref{cor-5} and 
Remark~\ref{r-2} together with the fact that for any 
$\Gamma\subset \bN$ such that each of the sequences 
$\Big\{\mu_{_{V_{p,n(p)}}}\Big\}_{p\in\Gamma}$, 
$\Big\{\mu_{_{S_{p,n(p)}}}\Big\}_{p\in\Gamma}$
and $\Big\{\mu_{_{Q_{2,p}}}\Big\}_{p\in\Gamma}$ weak-$^*$ converges the 
following holds:
\begin{equation*}
\begin{split}
&\lim_{\Gamma\ni p\to\infty}
\frac{(p-k)\al_k(n,p,k)}{(p-1)\al_k(n,p,k)+n(p)-1}=1,\\
&\frac{k-1}{p-1}\le 
\frac{(p-k)\al_k(n,p,k)+n(p)-1}{(p-1)\al_k(n,p,k)+n(p)-1}\le
\frac{k}{p}\,\text{ if }p\ge k,\\
&\left|\tilde{\mu}_{_{Q_{2,p}}}(f)-\mu_{_{Q_{2,p}}}(f)\right|\le 
\sum_{l=1}^{p}\frac{|\al_k(n,p,k)-pa_l-n(p)+1|}{p[(p-1)\al_k(n,p,k)+n(p)-1]}
\max_{z\in K}|f(z)|\\
&\phantom{\left|\tilde{\mu}_{_{Q_{2,p}}}(f)-\mu_{_{Q_{2,p}}}(f)\right|}\,
\le \frac{2}{p-1}\max_{z\in K}|f(z)|,
\end{split}
\end{equation*}
where $K$ is a compact subset of $\bC$ as in 
Corollary~\ref{cor-7}, $f:K\to \bR$ is any (continuous) convex function,
$\mu_{_{Q_{2,p}}}$ denotes the root-counting measure of $Q_{2,p}$ and
$\tilde{\mu}_{_{Q_{2,p}}}\in\calM(K)$ is defined by
\begin{equation*}
\begin{split}
&\su\big(\tilde{\mu}_{_{Q_{2,p}}}\big)
=\calZ(Q_{2,p})=\{\ze_{l,p}\}_{l=1}^{p},\\
&\tilde{\mu}_{_{Q_{2,p}}}(\{\ze_{l,p}\})
=\frac{\al_k(n,p,k)-a_l}{(p-1)\al_k(n,p,k)+n(p)-1},\quad 1\le l\le p.
\end{split}
\end{equation*}

Turning to the proof of Theorem~\ref{t-j} recall first (cf., e.g., \cite{Sz})
that $P_n^{(\al,\be)}$ satisfies Jacobi's equation, i.e., the homogeneous
second order linear differential equation
$$(1-z^2)y''(z)+[\be-\al-(\al+\be+2)z]y'(z)+n(n+\al+\be+1)y(z)=0.$$
Therefore, if $\calZ\Big(P_n^{(\al,\be)'}\Big)=\{\ze_{n,j}'\}_{j=1}^{n-1}$
denotes the zero set of $P_n^{(\al,\be)'}$ then
\begin{multline*}
-\frac{n(n+\al+\be+1)P_n^{(\al,\be)}(z)}{(1-z^2)P_n^{(\al,\be)'}(z)}
=\frac{P_n^{(\al,\be)''}(z)}{P_n^{(\al,\be)'}(z)}
+\frac{\be-\al-(\al+\be+2)z}{1-z^2}\\
=\sum_{j=1}^{n-1}\frac{1}{z-\ze_{n,j}'}+\frac{\al+1}{z-1}+\frac{\be+1}{z+1}.
\end{multline*}
Since $\al>-1$ and $\be>-1$ we may apply Lemma~\ref{l-1} to get
\begin{equation}\label{eq-j-strong}
(n+\al+\be+1)\sum_{i=1}^{n}f(\ze_{n,i})\le 
(n+\al+\be)\sum_{j=1}^{n-1}f\!\left(\ze_{n,j}'\right)+(n+\be)f(1)+(n+\al)f(-1)
\end{equation}
for any (continuous) convex function $f:[-1,1]\to \bR$ (cf.~Remark~\ref{r-1} in
\S \ref{s2}). Now by Corollary~\ref{c-cons} we know that
$$n\sum_{j=1}^{n-1}f\!\left(\ze_{n,j}'\right)\le 
(n-1)\sum_{i=1}^{n}f(\ze_{n,i}),$$
which together with~\eqref{eq-j-strong} immediately gives the desired 
conclusion.

\section{Further Results and Related Problems}\label{s4}

\subsection*{1.}\label{ss41}

One can actually obtain stronger albeit somewhat less transparent 
versions of Theorems~\ref{t-1}--\ref{t-2} and establish convex domination 
relations for vectors whose coordinates are symmetric functions
on (subsets of) $\calZ(V)\vee \calZ(S)$ and 
$\calZ(S^{(k-1)})\vee \calZ(Q_2)$, respectively. Indeed, given $d\in\bN$
and $e\in\bZ_+$ with $e\le d$ let $\Pi_{d,e}$ denote the $e$-th 
elementary symmetric function on $d$ elements. Then \cite[Corollary 3]{Bo} 
shows that the zeros and poles of the function $\Phi$ defined in \eqref{phi} 
satisfy the following inequalities.

\begin{lemma}\label{l-2}
Let $\Phi$ be as in \eqref{phi} and denote its zeros by $w_j$, 
$1\le j\le m-1$, where as before it is understood that $z_i$ counts as 
a ``zero'' of $\Phi$ of multiplicity $m_i-1$ if it occurs precisely 
$m_i$ times in \eqref{phi}. Then
\begin{multline*}
\sum_{1\le j_1<\ldots<j_d\le m-1}
f\big(\Pi_{d,e}(w_{j_1},\ldots,w_{j_d})\big)\\
\le \sum_{1\le i_1<\ldots<i_d\le m}\!\left(1-\sum_{l=1}^{d}\tau_{i_l}\right)
f\big(\Pi_{d,e}(z_{i_1},\ldots,z_{i_d})\big)
\end{multline*}
for any convex function $f:\bC\to \bR$, $d\in \{1,\ldots,m-1\}$ and
$e\in\bZ_+$, $e\le d$.
\end{lemma}

Note that Lemma~\ref{l-1} corresponds to the case when $e=d=1$ in 
Lemma~\ref{l-2}. Using the latter with $m=n-k+p+1$ and 
\begin{equation*}
\begin{split}
&\tau_j=
\begin{cases}
\al(n,p,k)^{-1},\quad 1\le j\le n-k+1,\\
\al(n,p,k)^{-1}a_{j-n+k-1},\quad n-k+2\le j\le n-k+p+1,
\end{cases}\\
&\{z_i\}_{i=1}^{n-k+p+1}=\calZ(S^{(k-1)})\vee \calZ(Q_2),\quad
\{w_j\}_{j=1}^{n-k+p}=\calZ(V)\vee \calZ(S)
\end{split}
\end{equation*}
together with \eqref{gen-k} one can then deduce weighted majorization relations
of the aforementioned type that strengthen \eqref{eq-strong} in various 
ways. For special choices of the function $f$ one can slightly simplify these 
relations by first separating elements in $\calZ(S^{(k-1)})$ from those in 
$\calZ(Q_2)$ (once the terms occurring in the right-hand side of the 
above inequality are appropriately regrouped) and then using 
Corollary~\ref{c-cons} in order to 
estimate from above all resulting expressions that contain elements in 
$\calZ(S^{(k-1)})$ by means of similar expressions involving only elements
in $\calZ(S)$. Such simplifications can be made e.g.~for multiplicative 
convex functions of the form $f(z)=|z|^q$, $q\in\bZ_+$.

\subsection*{2.}\label{ss42}

In view of the above results it would be interesting to know whether
similar properties with respect to the Choquet order also hold
for spectral polynomials of more general classes of (Lam\'e-like) operators. 
Let 
\begin{equation}\label{hlo}
\dd(z)=\sum_{i=m}^{k}Q_i(z)\frac{d^i}{dz^i}
\end{equation}
be a linear ordinary differential operator of order $k$ with polynomial 
coefficients.
Following the terminology that we already employed for \eqref{h-d} 
(cf.~\cite{BBS}) we call
$\dd(z)$ a {\em higher Lam\'e operator} if its {\em Fuchs index}
$r:=\max_{m\le i\le k}(\deg Q_i-i)$ is non-negative. If $r=0$ then $\dd(z)$
is usually referred to as an {\em exactly solvable operator} in the physics 
literature. A higher Lam\'e operator $\dd(z)$ given by \eqref{hlo} is said to 
be {\em non-degenerate} if $\deg Q_k=k+r$, which is equivalent to the (quite
natural) requirement that $\dd(z)$ has either a regular or a regular singular 
point at $\infty$. For such an operator one may then consider the 
multiparameter spectral problem stated in \eqref{k-lame} and the corresponding
notions of ($n$-)solvability and higher Lam\'e (i.e., Van Vleck and 
Heine-Stieltjes) polynomials. A systematic study of the latter was recently 
made in \cite{BBS}. In particular, in {\em op.~cit.~}it was shown that
a non-degenerate higher Lam\'e operator $\dd(z)$ is $n$-solvable for all
sufficiently large $n$ and it was further proved that if the coefficients of 
$\dd(z)$ are algebraically independent then 
for any $n\in\bN$ there are exactly $\binom{n+r}{n}$ Van Vleck
polynomials and as many degree $n$ Heine-Stieltjes
polynomials, thus generalizing Heine's result 
(cf.~Remarks~\ref{r-appl-2} and~\ref{r-appl-k}). 

\begin{problem}\label{pb1}
Extend Theorems~\ref{t-1} and~\ref{t-2} to non-degenerate higher Lam\'e 
operators (subject to appropriate conditions).
\end{problem}

An important class of linear operators which seems particularly 
well-suited for Problem~\ref{pb1} consists of non-degenerate higher Lam\'e 
operators that also preserve hyperbolicity (HPOs). Indeed, as we already 
mentioned in Remark~\ref{r-appl-k} a complete classification of all HPOs
-- i.e., linear operators $T$ on $\bR[z]$ such that $T(P(z))$ has all 
real zeros whenever $P\in\bR[z]$ has all real zeros -- was recently obtained
in \cite{BBS1}. Moreover, various properties and characterizations of 
HPOs that belong to the Weyl algebra $\calA_1$ (that is, operators of the
form~\eqref{hlo}) were established in \cite{BBS2}. For instance, in
{\em op.~cit.~}it was shown that the coefficients $Q_i(z)$ of such an operator
have all real zeros and satisfy interlacing properties like those 
in~\eqref{not-1}. (Note e.g.~that if $Q_2(z)$, $Q_1(z)$ are as in~\eqref{not-1}
and $\calZ(Q_2)\subset \bR$ then the corresponding operator $\dd(z)$ given
by~\eqref{h-d} is an HPO.) Furthermore, in \cite{BBS} it was proved that
if a non-degenerate higher Lam\'e operator $\dd(z)$ is also an HPO then 
$\dd(z)$ is solvable and all its Van Vleck and
Heine-Stieltjes polynomials have simple real zeros. Finally, 
\cite[Conjecture 1]{Bo1} claims that HPOs either preserve or reverse 
the Choquet order on real-zero polynomials and that in particular, 
if $\dd(z)$ is 
an HPO of the form \eqref{hlo} then 
$\calZ\big(\dd(z)P(z)\big)\pr \calZ\big(\dd(z)Q(z)\big)$
whenever $P,Q\in\bR[z]$, $\deg P=\deg Q$, 
$\calZ(P),\calZ(Q)\subset\bR$, $\calZ(P)\pr\calZ(Q)$. 
For results supporting this conjecture, see \cite{Bo1,Bo2,BS0}. 
Problem~\ref{pb1} should therefore be 
particularly interesting 
for non-degenerate higher Lam\'e operators of HPO type.

\subsection*{3.}\label{ss43}

Let $P_n(z)$, $n\in\bZ_+$, be polynomials orthogonal with respect to a 
weight function $\om$ supported on a (finite or infinite) interval $[a,b]$
with $\om(z)>0$, $z\in (a,b)$. It is well known that such polynomial 
families satisfy a $3$-term recurrence relation 
$$zP_n(z)=a_{n+1}P_{n+1}(z)+b_nP_n(z)+a_nP_{n-1}(z)$$
and that under some mild assumptions on the weight function $\om$
(see, e.g., \cite{I}) they also satisfy a differential recurrence relation
of the type
$$P_n'(z)=A_n(z)P_{n-1}(z)-B_n(z)P_n(z)$$ 
and a second order differential equation of the form
$$P_n''(z)+R_n(z)P_n'(z)+S_n(z)P_n(z)=0,$$
where $\{A_n(z)\}_{n\in\bZ_+}$, $\{B_n(z)\}_{n\in\bZ_+}$, 
$\{R_n(z)\}_{n\in\bZ_+}$, $\{S_n(z)\}_{n\in\bZ_+}$ 
are certain function sequences. It is therefore natural to ask the following 
questions.

\begin{problem}\label{pb2}
Investigate whether there are weighted majorization relations similar to 
those in Theorem~\ref{t-j} and$\,/\,$or weighted majorization relations 
involving the zeros of any two (or three) consecutive terms for 
\begin{itemize}
\item[(i)] Laguerre-like and Hermite-like polynomials; 
\item[(ii)] (appropriate classes of) general orthogonal polynomials.
\end{itemize}
\end{problem}

\end{document}